\documentclass[12pt]{article}
\usepackage{amssymb}
\usepackage{latexsym,bm}
\usepackage{graphicx}
\usepackage{amsmath}

\newcommand{\bea}{\begin{eqnarray*}}
\newcommand{\eea}{\end{eqnarray*}}
\newcommand{\be}{\begin{equation}}
\newcommand{\ee}{\end{equation}}
\newcommand{\ben}{\begin{eqnarray*}}
\newcommand{\een}{\end{eqnarray*}}

\date{}
\begin{document}
\title{Detour-saturated graphs of small girths\footnote{E-mail addresses:
{\tt 235711gm@sina.com}(P.Qiao),
{\tt zhan@math.ecnu.edu.cn}(X.Zhan).}}
\author{Pu Qiao, Xingzhi Zhan\thanks{Corresponding author.}\\
{\small Department of Mathematics, East China Normal University, Shanghai 200241, China}
 } \maketitle
\begin{abstract}
 A detour of a graph G is a longest path in G. The detour order of G is the number of vertices in a detour of G.
 A graph is said to be detour-saturated if the addition of any edge increases strictly the detour order. L.W. Beineke,
 J.E. Dunbar and M. Frick asked the following three questions in 2005. (1) What is the smallest order of a detour-saturated graph
 of girth 4? (2) Let Pr be the graph obtained from the Petersen graph by splitting one of its vertices into three leaves.
 Is Pr the smallest triangle-free detour-saturated graph? (3) Does there exist a detour-saturated graph with finite girth bigger
 than 5? We answer these questions.
\end{abstract}


\section{Introduction}

We consider finite simple graphs. Denote by $V(G)$ and $E(G)$ the vertex set and the edge set of a graph $G$ respectively.
The cardinalities of $V(G)$ and $E(G)$ are called respectively the {\it order} and {\it size} of $G.$ The complement of $G$
is denoted by $\bar{G}.$

{\bf Definition 1.} A {\it detour} of a graph $G$ is a longest path in $G.$ The {\it detour order} of $G,$ denoted $\tau(G),$ is the number
of vertices in a detour of $G.$ $G$ is said to be {\it detour-saturated} if $\tau(G+e) > \tau(G)$ for any $e\in E(\bar{G}).$

For terminology we follow the pioneering work in [6], [7] and [1] on this topic. Clearly every graph is a spanning subgraph of
some detour-saturated graph with the same detour order. One motivation for studying detour-saturated graphs in [1] is to investigate
the famous path partition conjecture (PPC). It is observed in [1] that it is sufficient to prove PPC for detour-saturated graphs.

Given a graph $G$, let $x$ be a vertex with neighbors $y_1,y_2,\ldots,y_k.$ The graph $G$ {\it split} $x,$ denoted by $G_{s}[x],$
is the graph obtained from $G$ by replacing $x$ with $k$ independent vertices $x_1,x_2,\ldots,x_k$ and joining $x_i$ to $y_i,$
$i=1,2,\ldots,k.$ Let $Pr$ be the graph obtained from the Petersen graph by splitting any vertex. Given a class $\Omega$ of graphs,
a {\it smallest} graph in $\Omega$ is a graph of minimum order in $\Omega,$ having the minimum size among the graphs of that order in $\Omega.$

Among other results, Beineke, Dunbar and Frick [1, Theorem 3.2] proved that $Pr$ is the smallest detour-saturated graph with girth $5,$
and they [1, p.123] asked the following three questions:
\newline\indent (1) What is the smallest order of a detour-saturated graph of girth 4?
\newline\indent (2) Is $Pr$ the smallest triangle-free detour-saturated graph?
\newline\indent (3) Does there exist a detour-saturated graph with finite girth bigger than 5?

The purpose of this paper is to answer these questions.

\section{Main Results}

We have obtained the first two results (Theorems 1 and 2) by computer search using the software SageMath. In writing the programs
we also used the result in [1, Corollary 2.3] that if a detour-saturated graph has a vertex of degree $2,$ then it has a triangle.
We ran separate programs to verify that all the five graphs in Figures 1-4 below have the required properties.
It was observed in [1, p.123] that there exists a detour-saturated graph with girth $4$ and order $34$ by a result of Thomassen [8].

{\bf Theorem 1.} {\it The smallest order of a detour-saturated graph of girth $4$ is $14.$ There is a unique detour-saturated graph of
girth $4$ and order $14,$ which is depicted in Figure 1. It has detour order $13.$}

\vskip 3mm
\par
 \centerline{\includegraphics[width=4.0in]{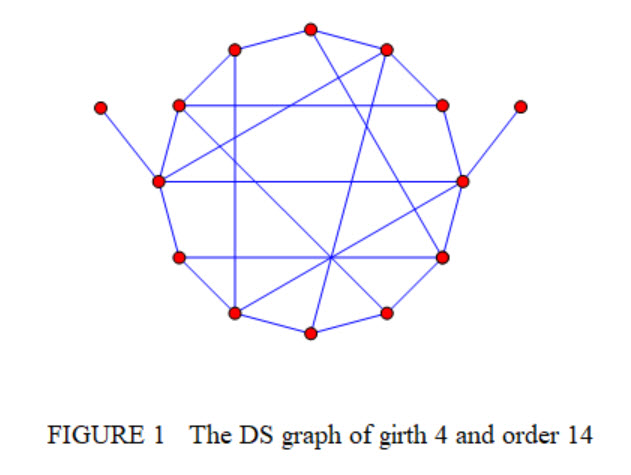}}
\par

{\bf Theorem 2.} {\it $Pr$ is the smallest triangle-free detour-saturated graph. There are exactly two triangle-free detour-saturated
graphs of order at most $12,$ which both have girth $5$ and order $12,$ have sizes $15$ and $17$ respectively, and are depicted
in Figure 2.}

\vskip 3mm
\par
 \centerline{\includegraphics[width=4.8in]{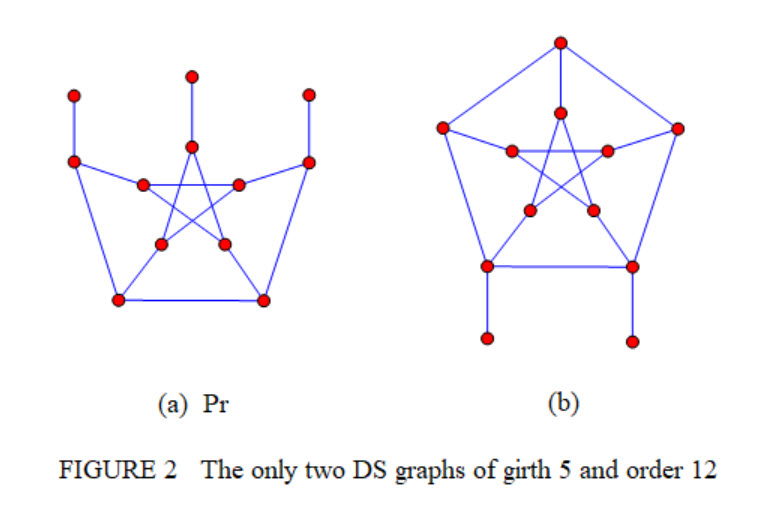}}
\par

Note that Theorems 1 and 2 above answer Questions 1 and 2 respectively and a little more. Theorem 4 below will answer Question 3.

{\bf Definition 2.} Let $G$ be a nonhamiltonian graph. $G$ is said to be {\it maximally nonhamiltonian} if $G+e$ is hamiltonian
for every edge $e\in E(\bar{G}).$ $G$ is {\it hypohamiltonian} if $G-v$ is hamiltonian for every vertex $v\in V(G).$
$G$ is {\it maximal hypohamiltonian} if $G$ is both maximally nonhamiltonian and hypohamiltonian.

We will need the following result in [1, Theorem 3.1].

{\bf Lemma 3} (Beineke, Dunbar and Frick) {\it Let $G$ be a maximal hypohamiltonian graph and let $x$ be a vertex of degree
$3$ in $G.$ Then $G_s[x]$ is a detour-saturated graph, with detour order equal to the order of $G.$}

{\bf Theorem 4.} {\it There exists an infinite family of detour-saturated graphs of girth $6$, and there exists a detour-saturated
graph of girth $7$ and order $30.$}

{\bf Proof.} Let $J_k$ be the cubic graphs of Isaacs [5] where $k\ge 3$ are odd integers. The definition of this class of graphs
can also be found in [3, pp.61-62]. Clark and Entringer [3, Corollary 11 and p.68] proved that  the graphs $J_k$ are
maximally nonhamiltonian and for $k\ge 5,$  $J_k$ are hypohamiltonian. Thus, for $k\ge 5,$ each $J_k$ is maximal hypohamiltonian.
It is known [3, p.68] that for $k\ge 7,$ the graphs $J_k$ have girth $6.$ From the construction of $J_k,$ it is clear that
when $k\ge 7,$ $J_k$ has at least two vertex-disjoint cycles of length $6,$ and consequently for any vertex $x$ of $J_k,$
the graph $J_k-x$ still has girth $6.$ By Lemma 3, we deduce that for every odd integer $k\ge 7$ and any vertex $x$ of $J_k,$
the graph $(J_k)_s[x]$ is a detour-saturated graph of girth $6.$ This proves the first statement of Theorem 4.

$(J_7)_s[x]$ has order $30$ and detour order $28.$ For one choice of the vertex $x,$ the graph $(J_7)_s[x]$ is depicted
in Figure 3.

\vskip 3mm
\par
 \centerline{\includegraphics[width=4.2in]{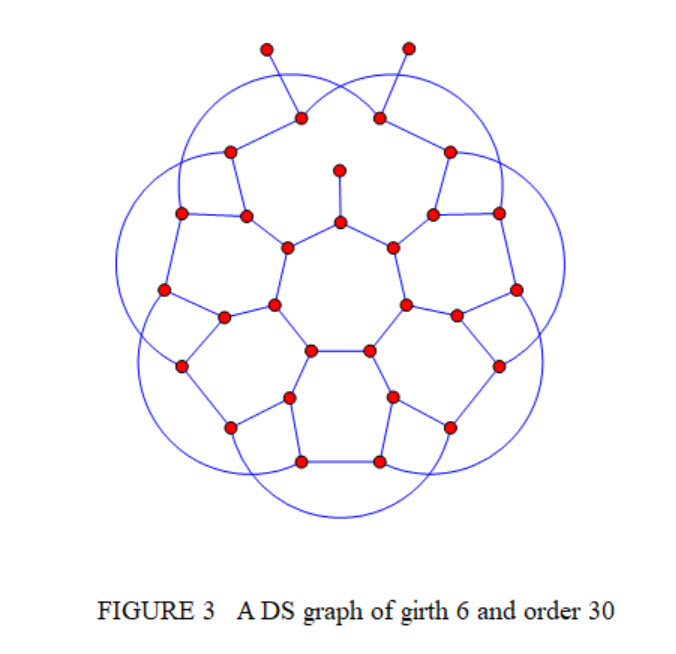}}
\par

Let $CT$ be the Coxeter graph whose definition can be found in [9] or [3, p.59]. It is a cubic graph of girth $7$ and order $28.$
Tutte [9] proved that $CT$ is nonhamiltonian. Bondy [2, pp.59-60] proved that $CT$ is hypohamiltonian and later Clark and Entringer
[3, Proposition 6] proved that $CT$ is maximally nonhamiltonian. Hence $CT$ is maximal hypohamiltonian.
Let $x$ be any vertex of $CT.$ Clearly $CT$ has at least two vertex-disjoint cycles of length $7,$ implying that $CT-x$ has girth $7.$
By Lemma 3, $CT_s[x]$ is a detour-saturated graph. $CT_s[x]$ has girth $7,$ order $30$ and detour order $28.$ For one choice of $x,$
$CT_s[x]$ is depicted in Figure 4.

\vskip 3mm
\par
 \centerline{\includegraphics[width=4.2in]{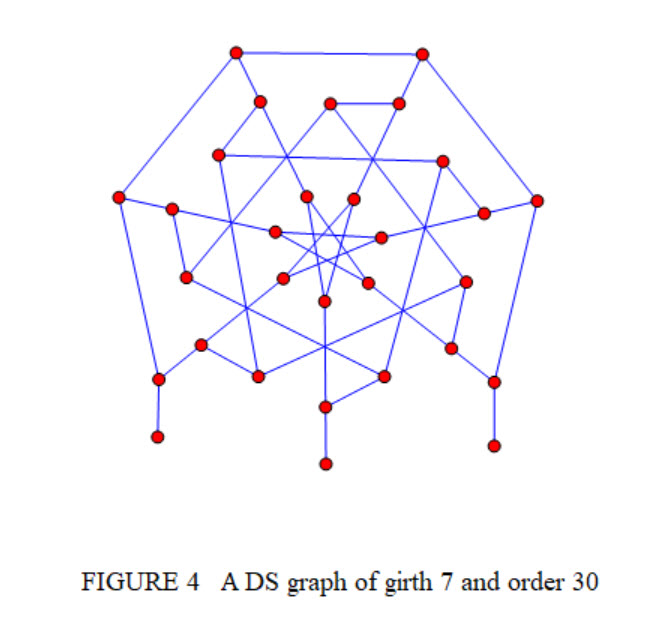}}
\par

The drawing in Figure 4 is a modification of a picture on the second page after the preface of the book [4].\hskip 5mm $\Box$

Finally we pose the following

{\bf Question.} Is it true that for every integer $g\ge 3,$ there exists a detour-saturated graph of girth $g$?

By the results in [1] and Theorems 1 and 4 above, the answer to this question is yes for $g=3,4,5,6,7.$

\vskip 5mm
{\bf Acknowledgement.}  This research  was supported by Science and Technology Commission of Shanghai Municipality (STCSM) grant 13dz2260400 and
 the NSFC grant 11671148.

\end{document}